\documentclass[a4paper,11pt]{amsart}
\usepackage[english]{babel}
\usepackage{amssymb,amsmath,amsfonts}
\usepackage{enumerate}
\usepackage{amsbib}

\newtheorem{theorem}{Theorem}

\newtheorem{statement}{Statement}

\newtheorem{corollary}{Corollary}
\newtheorem{problem}{Problem}

\theoremstyle{definition}
\newtheorem{definition}{Definition}

\newcommand\ind{\operatorname{ind}}
\newcommand\Ind{\operatorname{Ind}}

\begin{document}

\title{Coincidence of the Dimensions\\ of First Countable Spaces\\ with a Countable Network}

\author{I. M. Leibo}

\email{imleibo@mail.ru}

\address{Moscow Center for Continuous Mathematical Education, Moscow, Russia}

\begin{abstract}
The coincidence of the  $\Ind$ and $\dim$ dimensions for first countable paracompact $\sigma$-spaces is proved.
This gives a positive answer to A.~V.~Arkhangel'skii's question of whether the dimensions $\ind X$, $\Ind X$, and
$\dim X$ are equal for first countable spaces with a countable network.
\end{abstract}


\maketitle

\markboth{Leibo}{Coincidence of Dimensions}

All topological spaces considered in this paper are assumed to be normal $T_1$-spaces and called simply spaces.

In this paper, we use the terminology of the papers~\cite{1} and~\cite{3}.

There exist three classical topological dimension functions, $\Ind X$, $\dim X$, and $\ind X$ (see, e.g.,
\cite{1}); relationships between them for various classes of spaces is of special interest in dimension theory.
L.~A.~Tumarkin and W.~Hurewicz proved the equalities $\ind X=\dim X=\Ind X$ for a space $X$ with a countable base
\cite[Ch.~4, Sec.~8, Subsec.~2]{1}. Since every base is a network, it is quite natural to ask whether the
dimensions of a space with a countable network coincide; this question was posed by
Arkhangel'skii in \cite[Sec.~5]{2}. In the general case, the answer is negative. A space $X$ with a countable
network for which $\dim X=1$ but  $\Ind X=\ind X=2$ was constructed by Charalambous \cite{4} (note that the
dimensions $\Ind$ and $\ind$ are equal for any such space). When stating his coincidence problem, Arkhangel'skii
did not exclude the possibility that a counterexample may exist. Therefore, he simultaneously stated the
dimension coincidence problem for spaces with a countable network which satisfy certain additional conditions,
the most interesting of which are the presence of a topological group structure and the first axiom of
countability (i.e., the existence of a countable neighborhood base at every point). The question of whether $\ind
G= \dim G= \Ind G$ for any topological group $G$ with a countable network still remains open (for more than 50
years now).

In this paper, we give a positive answer to Arkhangel'skii's question of whether $\Ind X=\dim X=\ind X$ for any
first countable space with a countable network.

We also obtain a more general positive result, namely, that $\dim X = \Ind X$ for any first countable paracompact
$\sigma$-space $X$ (recall that $X$ is a $\sigma$-space if $X$ has a $\sigma$-discrete network, i.e., a network
being a countable union of discrete families of subsets of~$X$; in particular, any space with a countable network
is a $\sigma$-space).

The  main results of this paper are Theorems 1 and~2 and Corollary~1. The complete proofs of all statements are
given in~\cite{preprint}.

The coincidence of the dimensions $\dim X$, $\Ind X$, and $\ind X$ for a first countable space $X$ with a
countable network is implied by the following theorem.

 \begin{theorem}
\label{t1}
Let $X$ be a first countable  paracompact $\sigma$-space. Then the following conditions are equivalent:
\begin{enumerate}
\item[\rm(a)]
$\dim X \leq n$\textup;
\item[\rm(b)]
$\Ind X \leq n$\textup;
\item[\rm(c)]	
$X= \bigcup_{i=1}^{n+1}{X_i}$, where each~$X_i$  is a $G_\delta$-set and $\dim {X_i}\leq0$ for $i=1,2,\dots,
n+1$\textup;
\item[\rm(d)]
$X$ is the image of a first countable paracompact $\sigma$-space $X_0$ with $\dim X_0=0$ under a closed
continuous map $f\colon X_0 \to X$ such that $|f^{-1}(x)|\leq n+1$ for all~$x\in X$.
\end{enumerate}
\end{theorem}

We consider the equivalence of these particular conditions (a)--(d) because they are equivalent for any
metrizable space.

It follows from the equivalence of conditions (a) and  (b) in Theorem~1 that $\dim X=\Ind X$
for any first countable  paracompact $\sigma$-space $X$. As is known \cite[Sec.~8, Lemma~2]{1}, for any space $X$
with a countable network, we have $\dim X\leq \ind X\leq \Ind X$. This gives the following corollary.

 \begin{corollary}
\label{c1}
Let $X$ be a first countable space with a countable network. Then $\ind X= \dim X=\Ind X$.
\end{corollary}

Theorem~1 follows from Theorem~2 given below. Before stating it, we give two definitions.

\begin{definition}[\cite{6}]
Given a normal space $Z$ of finite dimension $\Ind Z=n$,  we say that a closed subset $F$ of $Z$ and its open
neighborhood $OF$ \emph{determine the dimension} $\Ind Z$ of $Z$ if, for any open neighborhood $U$ of $F$ such that
$F \subseteq U \subseteq \overline{U}\subseteq OF$, we have $\Ind\operatorname{Fr}U \ge n-1$.
\end{definition}

\begin{definition}[\cite{6}]
A family $\varphi=\{F_{\alpha}, OF_{\alpha}\}_{\alpha \in A}$ of closed  subsets $F_{\alpha}$ of a normal space $Z$
and their open neighborhoods $OF_{\alpha}$ is called an \emph{everywhere f-system} (an \emph{f-system})
if the family $\{OF_{\alpha}\}_{\alpha \in A}$ is $\sigma$-discrete and, given any (any closed) set $M \subseteq
Z$, there exists an index $\alpha_0 \in A$ for which the pair $\{M \cap F_{\alpha_0}$, $M\cap OF_{\alpha_0}$
determines the dimension $\Ind M$.
\end{definition}

Note that if the space $Z$ in Definition~2 is a paracompact $\sigma$-space, then it suffices to require the family
$\{OF_{\alpha}\}_{{\alpha}{\in}{A}}$ to be $\sigma$-locally finite in~$Z$ ~\cite[Coroll.~1.11]{7}.

Using the notion of an f-system, we can state a general condition sufficient for the coincidence of the
dimensions $\dim$ and~$\Ind$.

\begin{statement}[{\cite[Theorem~1]{6}}]
If a normal space $X$ has an f-system, then $\Ind X=\dim X$.
\end{statement}

Thus, the existence of an f-system in a first countable paracompact $\sigma$-space $X$ implies
$\dim X=\Ind X$ (and the equivalence of conditions (a) and  (b) in Theorem~1). Therefore,
we have $\ind X = \dim X = \Ind X$ for a  first countable space $X$ with a countable network
(Corollary~1).

According to~\cite{6}, the existence  of an everywhere f-system in a first countable paracompact $\sigma$-space
$X$ is sufficient for the equivalence of conditions (a)--(d)  in Theorem~1.

Thus, to prove Theorem~1, it suffices to prove the existence of an everywhere f-system in any first countable
paracompact $\sigma$-space. It is straightforward to construct such a system in a metrizable space, but the
case of a first countable paracompact $\sigma$-space is much more complicated.

However, there exists a fairly large subclass of the class of paracompact $\sigma$-spaces in which everywhere
f-systems surely exist, and it contains all first countable paracompact $\sigma$-spaces. This is the class of
almost semicanonical spaces. It was essentially introduced by Dugundji \cite{5}; the term ``semicanonical'' was
suggested by Robert Cauty.

\begin{definition} Let $A$ be a closed subset of a space $X$. The pair $(X,A)$ is said to be \emph{semicanonical}
if there exists an open cover $\omega$ of the open set $X\setminus A$ satisfying the following condition: any
open neighborhood $OA$ of $A$ contains an open neighborhood $UA\subseteq OA$ of $A$ such that any
element $G$ of $\omega$ intersecting $UA$ is contained in $OA$.

A space $X$ in which the pair $(X,A)$ is  semicanonical for any closed $A\subset X$ is called a
\emph{semicanonical space}.
\end{definition}

Let $X$ be a paracompact $\sigma$-space, and let $\gamma=\{F_{\alpha}: \alpha\in A\}$ be a $\sigma$-discrete
network of $X$; suppose that $A=\bigcup_{i=1}^{\infty} A_{i}$ and the family
$\gamma_i=\{F_{\alpha}: \alpha \in A_{i}\}$ is
discrete in $X$ for each $i=1,2, \dots$\,. Since $X$ is regular, we can assume that all $F_\alpha$ are closed.
We set $F_i = \bigcup\{F_{\alpha}: \alpha \in A_{i}\}$ for $i=1,2, \dots$\,.

\begin{definition}
A paracompact $\sigma$-space $X$ is said to be \emph{almost semicanonical} if it has a network $\gamma$ of the
above form such that every pair $(X, F_i)$ is semicanonical.
\end{definition}

\begin{statement}[{\cite[Theorem 1.17]{7}}]
Any  almost semicanonical space $X$ contains an f-system and an everywhere f-system.
\end{statement}

To prove Theorem~1, it remains to apply the following theorem.

\begin{theorem}
Any first countable paracompact $\sigma$-space $X$ is almost semicanonical.
\end{theorem}

 The results stated above give a general scheme for proving the coincidence of the dimensions $\dim X$ and $\Ind
X$ for almost semicanonical spaces $X$, which constitute a fairly large class. This class contains all metrizable
spaces, Nagata (that is, stratifiable first countable) spaces and our first countable paracompact $\sigma$-spaces.

The class of almost semicanonical spaces does not contain all spaces with a countable network; e.g., it does not
contain the space $X$ with a countable network and noncoinciding dimensions that was constructed by Charalambous
\cite{4}. There also exist stratifiable spaces which are not semicanonical~\cite{8}.

In proving the equality $\dim X = \Ind X$ for a given space $X$ by the scheme described above,
most difficult is verifying the semicanonicity of $X$. The proof is very involved even for a first countable
paracompact $\sigma$-space $X$ satisfying the very strong additional assumption of the existence on~$X$ of a
semimetric continuous with respect to one variable~\cite[Theorem~2.15]{7}.

 Still, this is a working method for proving the equality $\dim X = \Ind X$, which suggests the following
problem.

 \begin{problem}
Let $G$ be a topological group with a countable network. Is it true that the underlying space of $G$ is
almost semicanonical?
\end{problem}

 A positive answer to this question would give a positive answer to Arkhan\-gel'skii's question of whether the
dimensions of a topological group with a countable network coincide.

 It is my pleasure to express gratitude to Professor A.~V.~Arkhangel'skii for attention and to Professor
 O.~V.~Sipacheva for useful discussions.

\def\refname{References}

\end{document}